\documentclass[A4paper,11pt,twoside,english]{article}
\usepackage{babel}
\usepackage{graphicx}
\usepackage{latexsym}
\usepackage{amsmath}
\usepackage{makeidx}
\usepackage[utf8]{inputenc}
\usepackage[utf8]{inputenc}
\usepackage{amsfonts}
\usepackage{amssymb}
\usepackage{amsmath}
\usepackage{amsthm}
\usepackage{tocbibind}
\usepackage{tcolorbox}
\usepackage[framemethod=TikZ]{mdframed}
\usepackage{lipsum}
\usepackage[mathscr]{eucal}
\usepackage{graphicx}
\usepackage{pdfpages}
\usepackage{color}
\usepackage{array}
\usepackage{hyperref}
\usepackage[rightcaption]{sidecap}
\usepackage[noadjust]{cite}

\title{\bf A characterization of generalized Lipschitz classes by the rate of convergence of semi-discrete operators
}
\author{ {\bf Danilo Costarelli}\thanks{Corresponding author} \hskip1cm {\bf Michele Piconi}\hskip1cm {\bf Gianluca Vinti}
\\
Department of Mathematics and Computer Science \\
University of Perugia\\
1, Via Vanvitelli, 06123 Perugia, Italy
\\
{\small {\tt danilo.costarelli@unipg.it}} - {\small {\tt michele.piconi@unipg.it}} - {\small {\tt gianluca.vinti@unipg.it}}
\date{}
 }

\newtheorem{lemma}{Lemma}[section]
\newtheorem{teorema}{Theorem}[section]
\newtheorem{cor}{Corollary}[section]
\theoremstyle{definition}
\newtheorem{remark}{Remark}[section]

\usepackage{fancyhdr}
\begin{document}
\maketitle
\begin{abstract}
In this paper, we establish a comprehensive characterization of the generalized Lipschitz classes through the study of the rate of convergence of a family of semi-discrete sampling operators, of Durrmeyer type, in $L^p$-setting. To achieve this goal, we provide direct approximation results, which lead to quantitative estimates based on suitable $K$-functionals in Sobolev spaces and, consequently, on higher-order moduli of smoothness. Additionally, we introduce a further approach employing the celebrated Hardy-Littlewood maximal inequality to weaken the assumptions required on the kernel functions. These direct theorems are essential for obtaining qualitative approximation results in suitable Lipschitz and generalized Lipschitz classes, as they also provide conditions for studying the rate of convergence when functions belonging to Sobolev spaces are considered. The converse implication is, in general, delicate, and actually consists in addressing an inverse approximation problem allowing to deduce regularity properties of a function from a given rate of convergence. Thus, through both direct and inverse results, we establish the desired characterization of the considered Lipschitz classes based on the $L^p$-convergence rate of Durrmeyer sampling operators. Finally, we provide remarkable applications of the theory, based on suitable combinations of kernels that satisfy the crucial Strang-Fix type condition used here allowing to both enhance the rate of convergence and to predict the signals.
\end{abstract}

\section{Introduction}
Let $f$ be a continuous and real-valued function on the closed unit interval, namely $f\in C[0,1]$, and let
\[
\left(B_nf\right)(x):=\sum_{\nu=0}^nf\left(\dfrac{\nu}{n}\right)p_{\nu,n}(x),\quad p_{\nu,n}(x)=\binom{n}{\nu}x^\nu(1-x)^{n-\nu},\quad x\in[0,1],
\]
be its associated Bernstein polynomial. It is well-known that Bernstein polynomials represent one of the first considered families of approximation operators (of the polynomial type), which are still widely studied by several researchers, along with their numerous and various generalizations (see \cite{butzer2009voronovskaya,Heilmann,karsli2021approximation}). In literature, several problems regarding such polynomials have been approached, consisting, other than convergence results, in both quantitative and qualitative estimates for the order of approximation in case of functions belonging to suitable Lipschitz classes. Among these, we recall the following very classical result 
\[
\left|\left(B_nf\right)(x)-f(x)\right|\le M\left(\frac{x(1-x)}{n}\right)^{\frac{\alpha}{2}},\quad 0<\alpha\le 1,\quad x\in[0,1],
\]
(see, e.g., \cite{Lorentz53,Lorentz63}). Such estimate shows that the order of approximation is significantly influenced by the regularity properties of $f$. Specifically, stronger regularity properties, such as Lipschitz conditions, are associated with higher orders of approximation. Furthermore, as is common in Approximation Theory, the properties of the involved kernels (in $B_nf$ are the $p_{\nu,n}$'s) are crucial in this context (see, e.g., \cite{FAA71}).
\\
On the other hand, inverse statements that investigate the structural regularity properties of a function from the knowledge of the order of approximation are much more difficult to prove. This problem has been solved for Bernstein polynomials in \cite{Berens72,BeckerNessel78,Becker79}.
\\
In order to face both direct and inverse problems in more general functional spaces, we first need to generalize the classical definition of Bernstein polynomials so that it makes sense for a larger class of functions, even not necessarily continuous. In this direction, the Durrmeyer's method turned out to be very satisfactory. Durrmeyer's idea was to replace the sample values $f\left(\dfrac{\nu}{n}\right)$ by a general convolution integral, where the same generating kernel $p_{\nu,n}$ appears, i.e.,
\[
\left(D_nf\right)(x):=(n+1)\sum_{\nu=0}^np_{\nu,n}(x)\int_{0}^1p_{\nu,n}(u)f(u)du,\quad x\in[0,1].
\]
The literature about these operators is considerably wide and in development: we quote, e.g., \cite{Derriennic, Berdysheva, Gonska, Heilmann}.\\ In recent years, this method was also applied to the generalized sampling series, giving basic mathematical models for image and signal processing. The theory of generalized sampling series has been firstly introduced by the Paul Butzer’s school in Aachen during the 1980s, opening the way for several families of sampling type operators \cite{17, BS93}. In particular, applying the Durrmeyer method to the generalized sampling series, we obtain semi-discrete sampling operators known in literature as \textit{Durrmeyer sampling type operators}, of the form
\begin{equation*}
\left(\mathcal{D}_w^{\varphi,\psi}f\right)(x):=\sum_{k\in\mathbb{Z}}\varphi(wx-k)w\int_{\mathbb{R}}\psi(wu-k)f(u)du,\quad w>0,\quad x\in\mathbb{R},
\end{equation*}
where $\varphi$ and $\psi$ are the generating kernels (discrete and continuous, respectively). Such operators have been first introduced in \cite{12} to establish an asymptotic formula for functions with a polynomial growth, obtaining, as a consequence, also a Voronovskaja type formula (see also \cite{HR,CMGR2014,2,bajpeyi2022approximation}). Furthermore, due to their general formulation, we point out that Durrmeyer sampling type operators include other remarkable families of sampling type operators. Among these, we recall the well-known \textit{generalized sampling type operators} (\cite{BS93, 17,ITSF08}), that arise from a distributional version of Durrmeyer series where $\psi$ is the Dirac delta distribution. Furthermore, taking as $\psi$ the characteristic function of the interval $[0,1]$, we reduce to the celebrated \textit{Kantorovich sampling type operators} (\cite{AHR,ACCSV,2022JFAA, ABR,CCV23, acar2024characterization,draganov2024characterization}). 
\\
In our earlier papers, we approach the problem of convergence for Durrmeyer sampling type operators in the general framework of Orlicz spaces (\cite{BS22,BS23}), both in one and in more variables (see, respectively, \cite{EstimatesDurr,DurrConv} and \cite{MultiDurr}), including the $L^p$-spaces as particular case. Recently, in \cite{DurrRegularization} we investigate the regularization properties of Durrmeyer sampling type operators in $L^p$-spaces: this study turns out to be preparatory to approach a comprehensive characterization of the well-known generalized Lipschitz classes of order $\alpha>0$, via the rate of convergence of Durrmeyer sampling operators in the $L^p$-setting, that is precisely the object of study of the present paper.\\
It is important to note that the generalized Lipschitz classes we investigated are substantially equivalent to Sobolev spaces when $\alpha=r\in\mathbb{N}$ and $1<p\le+\infty$. However, the equivalence slightly differs when $p=1$, involving the space of bounded variation functions. To obtain the above characterization, we start from direct approximation results. Here, we introduce a crucial moment-type condition involving the algebraic moments of the kernels $\varphi$ and $\psi$. We also employed the well-known expansion in Sobolev spaces with integral remainder to obtain finer quantitative estimates for Durrmeyer sampling operators, using $K$-functionals defined in such setting. Moreover, we also introduce an alternative approach based on the celebrated Hardy-Littlewood maximal inequality to broaden the applicability of these results. Thus, we employ the equivalence between $K$-functionals and moduli of smoothness of the function to achieve a full estimation of the rate of convergence of the operators in the $L^p$-setting. Also a qualitative analysis has been carried out, when functions belonging to Lipschitz and generalized Lipschitz classes are considered.
\\
The converse implication is highly delicate and consists in deducing the regularity properties of a function from the knowledge of the rate of convergence of Durrmeyer sampling operators with respect to the $L^p$-norm. To this aim we use a preliminary lemma regarding the algebraic moments of the distributional derivatives of the discrete kernel and we argue by induction on the rate of convergence. 
Therefore, through a synthesis of direct and inverse results, we arrive to establish a full characterization of the generalized Lipschitz classes of order $\alpha>0$ within the $L^p$-setting, with $1\le p\le +\infty$.

\section{Preliminaries}
We denote by $L^p(\mathbb{R})$, $1\le p<+\infty$, the space of all Lebesgue measurable functions $f:\mathbb{R}\to\mathbb{R}$ such that the usual norm
\begin{equation*}
\|f\|_p:=\left(\int_{\mathbb{R}}|f(u)|^pdu\right)^\frac{1}{p}
\end{equation*}
is finite. If $p=+\infty$, we formally consider the space $C(\mathbb{R})$ of all uniformly continuous and bounded functions endowed with the usual sup-norm $\|\cdot\|_\infty$. Further, we denote by $BV(\mathbb{R})$ the space of functions of bounded variation on $\mathbb{R}$ and by $AC(\mathbb{R})$ that one of all absolutely continuous functions on $\mathbb{R}$.
\\
One of the main framework we use is represented by Sobolev spaces $W^{r,p}(\mathbb{R})$, $r\in\mathbb{N}$, defined as the collection of functions in $L^p(\mathbb{R})$ for which the (distributional) derivatives $f^{(r-1)}\in AC(\mathbb{R})$ and $f^{(r)}\in L^p(\mathbb{R})$. If $p=+\infty$, we consider the space $C^r(\mathbb{R})\subset C(\mathbb{R})$ of all $r$-fold differentiable functions on $\mathbb{R}$, with $f^{(i)}\in C(\mathbb{R})$, $i=1,\dots,r$. We may also introduce the analogous space $W^{r}(BV)$, replacing $L^p(\mathbb{R})$ by $BV(\mathbb{R})$ in the previous definition (see \cite{DVL}, p. 35).
\\
A natural framework to obtain qualitative estimates for the order of approximation is given by the \textit{Lipschitz spaces}. The latters can be also defined in the more general case of $L^p$-functions. Here, the space $Lip(\alpha,L^p)$, with $0<\alpha\leq 1$ and $1\le p\le+\infty$, is defined as the set of all functions $f\in L^p(\mathbb{R})$ for which the \textit{$L^p$-first order modulus of smoothness} satisfies
\begin{equation}\label{omega1}
\omega(f,\delta)_p:=\sup_{|t|\le\delta}\|f(\cdot+t)-f(\cdot)\|_p=\mathcal{O}\left(\delta^\alpha\right),\qquad\delta\to 0^+.
\end{equation}
In order to extend the definition of $Lip(\alpha,L^p)$ for $\alpha>1$, we present two different approaches. Let $\alpha>0$ be such that $\alpha=r+\beta$, where $r=0,1,\dots$ is an integer and $0<\beta\le 1$. The space denoted by $W^r Lip(\beta, L^p)$ consists of all functions of $L^p(\mathbb{R})$ for which $f,\dots,f^{(r-1)}$ are locally absolutely continuous on $\mathbb{R}$ and satisfy $f^{(r)}\in Lip(\beta, L^p)$.
\\
On the other hand, we may define the \textit{generalized Lipschitz space} $Lip^*(\alpha,p)$, with $\alpha>0$ and $1\le p\le+\infty$, introducing in (\ref{omega1}) the finite differences of higher order for $f$. 
\\
Considering $r>\alpha$ as the smallest integer greater than $\alpha$, that is $\lfloor\alpha\rfloor+1$ (where $\lfloor\cdot\rfloor$ denotes the integer part), $Lip^*(\alpha,p)$ will be the space of all functions $f\in L^p(\mathbb{R})$ for which the \textit{modulus of smoothness of order $r$} satisfies \[\omega_r(f,\delta)_p:=\sup_{|t|\le\delta}\|\Delta^r_t(f,\cdot)\|_p=\mathcal{O}(\delta^\alpha),\qquad\delta\to 0^+,\]
where
\[
\Delta^r_t(f,x):=\sum_{j=0}^r\binom{r}{j}(-1)^{r-j}f(x+jt),\qquad x\in\mathbb{R}.
\]
Among the useful properties of $\omega_r(f,\delta)$, we recall that for a given positive integer $k$
\[
\omega_{r+k}(f,\delta)_p\le \delta^r\omega_k\left(f^{(r)},\delta\right)_p,\qquad \delta\ge0.
\]
This immediately implies that
\begin{equation}\label{inclusion}
Lip^*(\alpha,L^p)\supseteq Lip(\alpha,L^p).
\end{equation}
In particular, for $0<\alpha\le 1$, inclusion (\ref{inclusion}) is actually an equality. It is also well-known that if $\alpha>0$ is not an integer, we have
\[
Lip^*(\alpha,L^p)= Lip(\alpha,L^p),
\]
while if $\alpha=r\in\mathbb{N}$, we have
\[
Lip(r,L^p)=\begin{cases} W^{r,p}(\mathbb{R}), & \mbox{if }\mbox{ $1<p\le+\infty$} \\ W^{r-1}(BV)\cap L^1(\mathbb{R}), & \mbox{if }\mbox{ $p=1$.}
\end{cases}
\\
\subseteq Lip^*(r,L^p).
\]
In any case, there holds
\[
\omega_r(f,\delta)_p\le M\delta^r,\qquad \delta>0,
\]
for every $1\le p\le+\infty$ and $f\in Lip(r,L^p)$.
\\
A useful tool to achieve suitable quantitative estimates in terms of $\omega_r(f,\cdot)_p$, when $f\in L^p(\mathbb{R})$ with $1\le p\le+\infty$, is the well-known \textit{K-functional}, defined as
\[
\mathcal{K}\left(f,t;L^p,W^{r,p}\right):=\inf_{g\in W^{r,p}(\mathbb{R})}\left\{\|f-g\|_p+t\|g^{(r)}\|_p\right\}\qquad t>0.
\]
The relation between $\omega_r(f,\cdot)_p$ and $\mathcal{K}\left(f,t;L^p,W^{r,p}\right)$ is given by the following classical result established by H. Johnen in 1972 (see, e.g., \cite{J72}).
\begin{teorema}\label{Johnen} For $1\le p\le+\infty$ and $r=1,2,\dots$ there exist two positive constants $C_1$, $C_2$ (depending only on $r$) such that for any $f\in L^p(\mathbb{R})$, it turns out that
\[
C_1\omega_r(f,t)_p\le \mathcal{K}\left(f,t^r;L^p,W^{r,p}\right)\le C_2\omega_r(f,t)_p,\qquad t>0.
\]
\end{teorema}

\noindent In order to face both direct and inverse problems for a general family of semi-discrete sampling type operators, we first have to introduce two fundamental functions $\varphi, \psi:\mathbb{R}\to\mathbb{R}$ that will play the role of \textit{kernels} for such operators. Basically, we suppose that both functions belong to $L^1(\mathbb{R})$, $\varphi$ is bounded on $\mathbb{R}$ and $\psi$ is bounded in a neighbourhood of the origin, as it is usual for the kernel functions.
\\
We define the \textit{discrete} and \textit{continuous algebraic moments} of $\varphi$ and $\psi$ of order $\nu$, respectively, as follows
\[ m_{\nu}(\varphi,u):=\sum_{k\in\mathbb{Z}} \varphi(u-k)(k-u)^\nu,\quad u\in\mathbb{R}, \ \nu \in \mathbb{N}_0,\] \\ and  \[\widetilde {m}_\nu\left(\psi\right):= \int_{\mathbb{R}}u^\nu \psi(u)du,\quad \nu \in \mathbb{N}_0\] \\and the \textit{discrete} and \textit{continuous absolute moments} as \\ 
\begin{equation*}\label{momentoassolutod} M_{\nu}(\varphi):=\sup_{u\in\mathbb{R} }\sum_{k\in\mathbb{Z}}\left | \varphi(u-k)\right|\left|u-k\right |^\nu,\quad \nu>0 \end{equation*} \\and\\  \[\widetilde {M}_\nu(\psi):= \int_{\mathbb{R}}\left |u\right|^\nu\left | \psi(u)\right|du,\quad \nu>0\]
respectively. In particular, if the following basic condition holds \[m_0(\varphi,u)=\widetilde{m}_0(\psi)=1,\] we call the functions $\varphi$ and $\psi$ as \textit{discrete} and \textit{continuous kernel}, respectively. 
\\
For $w>0$ and for kernels $\varphi$ and $\psi$, we can recall the definition of a family of operators by
 \begin{equation}\label{Durrmeyeroperator} \left ( \mathcal{D}_w^{\varphi,\psi}f\right)(x)= \sum_{k\in\mathbb{Z}} \varphi(wx-k) w\int_{\mathbb{R}}\psi(wu-k)f(u)du,\text{$\quad$} x\in \mathbb{R},
 \end{equation}
 for any given function $f:\mathbb{R}\to\mathbb{R}$ such that the above series is convergent on $\mathbb{R}$.
 $\mathcal{D}_w^{\varphi,\psi}$ are called the \textit{Durrmeyer sampling type operators} based upon $\varphi$ and $\psi$.\\
In order to achieve the main results herein presented, we additionally suppose the following conditions. First of all, we say that a kernel $\xi:\mathbb{R}\to\mathbb{R}$ satisfies condition $(\Theta)$, if there exists $\theta>0$ such that
\begin{equation}
\xi(u)=\mathcal{O}(|u|^{-\theta}),\quad|u|\to+\infty.\tag{$\Theta$}
\end{equation}
Let $r>1$ be an integer, we assume that
\begin{description}
\item[{(i)}] for every $\nu=1,\dots,r$, the algebraic moments $m_\nu(\varphi,u)$ with $\nu=1,\dots,r$ are independent of $u\in\mathbb{R}$. In this case, we simply write $m_\nu(\varphi,u)=m_\nu(\varphi)$;
\item[{(ii)}] for every $i=1,\dots,r-1$, the following sum vanishes
\[\sum_{\nu=0}^i\binom{i}{\nu}m_{i-\nu}(\varphi)\widetilde{m}_\nu(\psi)=0.\]
\end{description}
In particular, condition (ii) is actually a \textit{vanishing moment condition of Durrmeyer type} (also known as a \textit{Strang-Fix type condition}) that plays a crucial role in the present study. Later, we provide some examples of kernels for which the above assumptions (i) and (ii) are satisfied, showing that these conditions are not restrictive.
\begin{remark}\label{preliminaryremark}
We outline that if a function $\xi:\mathbb{R}\to\mathbb{R}$ satisfies $\xi(u)=\mathcal{O}(|u|^{-\theta})$, as $|u|\to+\infty$, with $\theta>\beta+1$, $\beta>0$, then
\[
M_\gamma(\xi)<+\infty,
\]
for every $0\le \gamma<\beta$. Moreover, if $\xi:\mathbb{R}\to\mathbb{R}$ is bounded on $\mathbb{R}$, then, for every $\gamma,\beta\in\mathbb{R}$ such that $0\le \gamma<\beta$, $M_{\beta}(\xi)<+\infty$ implies $M_{\gamma}(\xi)<+\infty$. In fact, it is easy to see that $M_\gamma(\xi)\le 2\|\xi\|_\infty+M_\beta(\xi)$. Clearly, if $\xi$ has compact support then $M_{\beta}(\xi)<+\infty$ for every $\beta\ge0$. In particular, the basic moment condition $M_0(\xi)<+\infty$ is easily ensured whenever $\xi$ satisfies ($\Theta$) with $\theta>1$.
\end{remark}
We now observe that if the discrete kernel $\varphi$ satisfies $M_0(\varphi)<+\infty$, then the Durrmeyer sampling type operators are well-defined, linear and bounded operators on $L^\infty(\mathbb{R})$.\\
On the other hand, in the $L^p$-setting, a similar situation holds. Indeed, we can state the following.
\begin{teorema}[\cite{DurrConv}]\label{thm2}
Let $\varphi$ and $\psi$ the discrete and continuous kernel respectively. 
 If $M_0(\varphi)+M_0(\psi)<+\infty$ and $1\le p\le+\infty$, then the operators
\[
\mathcal{D}^{\varphi,\psi}_w:L^p(\mathbb{R})\to L^p(\mathbb{R})
\]
are well-defined, linear and continuous. In particular, for every $f\in L^p(\mathbb{R})$, there holds
\[
\|\mathcal{D}^{\varphi,\psi}_wf\|_p\le\mu \|f\|_p,
\]
with
\[
\mu:=\begin{cases} M_0(\varphi)^{\frac{p-1}{p}}\|\psi\|_1^{\frac{p-1}{p}}\|\varphi\|_1^{\frac{1}{p}}M_0(\psi)^{\frac{1}{p}}, & \mbox{if }\mbox{ $1\le p<+\infty$} \\ M_0(\varphi)\|\psi\|_1, & \mbox{if }\mbox{ $p=+\infty$}.
\end{cases}
\]
Moreover, for every $f\in L^p(\mathbb{R})$, we have
\[
\lim_{w\to+\infty}\|\mathcal{D}^{\varphi,\psi}_wf-f\|_p=0.
\]
\end{teorema}

\section{Direct approximation results}\label{directsection}
Let $x\in\mathbb{R}$, by using the Taylor formula with integral remainder
\begin{equation*}
f(u)=f(x)+\sum_{i=1}^{r-1}\frac{f^{(i)}(x)}{i!}(u-x)^i+\int_x^u\frac{f^{(r)}(t)}{(r-1)!}(u-t)^{r-1}dt,\quad u\in\mathbb{R},
\end{equation*}
for every $f\in W^{r,p}(\mathbb{R})$ ($r>1$) if $1\le p<+\infty$ (see \cite{DVL}), or $f\in C^r(\mathbb{R})$ if $p=+\infty$, we obtain the following expansion for the Durrmeyer sampling type series
\begin{equation*}
\begin{split}
\left(\mathcal{D}_w^{\varphi,\psi}f\right)(x)&=f(x)+\sum_{i=1}^{r-1}\frac{f^{(i)}(x)}{i!}\sum_{k\in\mathbb{Z}}\varphi(wx-k)w\int_{\mathbb{R}}\psi(wu-k)(u-x)^idu
\\
&\qquad+\sum_{k\in\mathbb{Z}}\varphi(wx-k)w\int_\mathbb{R}\psi(wu-k)\left\{\int_x^u\frac{f^{(r)}(t)}{(r-1)!}(u-t)^{r-1}dt\right\}du, \quad x\in\mathbb{R}.
\end{split}
\end{equation*}
Let $x\in\mathbb{R}$ and $i=1,\dots,r-1$, we may write
\begin{equation*}
\begin{split}
I_i&:=\frac{f^{(i)}(x)}{i!}\sum_{k\in\mathbb{Z}}\varphi(wx-k)w\int_{\mathbb{R}}\psi(wu-k)(u-x)^idu
\\
&=\frac{f^{(i)}(x)}{i!w^i}\sum_{\nu=0}^i\binom{i}{\nu}\sum_{k\in\mathbb{Z}}\varphi(wx-k)(k-wx)^{i-\nu}w\int_{\mathbb{R}}\psi(wu-k)(wu-k)^\nu du
\\
&=\frac{f^{(i)}(x)}{i!w^i}\sum_{\nu=0}^i\binom{i}{\nu}m_{i-\nu}(\varphi)\widetilde{m}_\nu(\psi).
\end{split}
\end{equation*}
Therefore, using condition (ii), we have
\begin{align}\label{taylor}
\left(\mathcal{D}_w^{\varphi,\psi}f\right)(x)
&=f(x)+\sum_{i=1}^{r-1}I_i+\sum_{k\in\mathbb{Z}}\varphi(wx-k)w\int_\mathbb{R}\psi(wu-k)
\left\{\int_x^u\frac{f^{(r)}(t)}{(r-1)!}(u-t)^{r-1}dt\right\}du\nonumber
\\
&=f(x)+\sum_{i=1}^{r-1}\frac{f^{(i)}(x)}{i!w^i}\sum_{\nu=0}^i\binom{i}{\nu}m_{i-\nu}(\varphi)\widetilde{m}_\nu(\psi)\nonumber
\\
&\qquad+\sum_{k\in\mathbb{Z}}\varphi(wx-k)w\int_\mathbb{R}\psi(wu-k)\left\{\int_x^u\frac{f^{(r)}(t)}{(r-1)!}(u-t)^{r-1}dt\right\}du\nonumber
\\
&=f(x)+\sum_{k\in\mathbb{Z}}\varphi(wx-k)w\int_\mathbb{R}\psi(wu-k)\left\{\int_x^u\frac{f^{(r)}(t)}{(r-1)!}(u-t)^{r-1}dt\right\}du,
\end{align}
$x\in\mathbb{R}$. This equality (\ref{taylor}) will be fundamental in the following.
\begin{remark}\label{rmk1} Let $f\in L^p(\mathbb{R})$ and $g\in W^{r,p}(\mathbb{R})$, $1\le p\le+\infty$. Assuming that $M_0(\varphi)+M_0(\psi)<+\infty$, there exists a suitable constant $\mu>0$ such that
\begin{equation}
\begin{split}
\|\mathcal{D}_w^{\varphi,\psi}f-f\|_p&\le\|\mathcal{D}_w^{\varphi,\psi}f-\mathcal{D}_w^{\varphi,\psi}g\|_p+ \|\mathcal{D}_w^{\varphi,\psi}g-g\|_p+\|g-f\|_p
\\
&\le(\mu+1)\|f-g\|_p+\|\mathcal{D}_w^{\varphi,\psi}g-g\|_p,
\end{split}
\end{equation}
where $\mu>0$ is the constant given in Theorem \ref{thm2}.
\end{remark}
\noindent We stress that, from now on, $r\in\mathbb{N}$ always denotes an integer strictly greater than 1. Now, we are ready to prove what follows.
\begin{teorema}\label{Linf} Let $\varphi$ and $\psi$ be the discrete and continuous kernel respectively satisfying $(i)$, $(ii)$ and condition $(\Theta)$ with $\theta>r+1$. Then, for any $f\in C(\mathbb{R})$ there exists a constant $\Lambda_\infty>0$ such that
\[
\|\mathcal{D}_w^{\varphi,\psi}f-f\|_\infty\le \Lambda_\infty\cdot\omega_r\left(f,\frac{1}{w}\right)_\infty,\qquad w>0.
\]
\end{teorema}
\begin{proof}
Let $f\in C(\mathbb{R})$ be fixed and $g\in C^{r}(\mathbb{R})$. By using the Taylor expansion with integral remainder for $g$ as in (\ref{taylor}), we obtain
\[
\left(\mathcal{D}_w^{\varphi,\psi}g\right)(x)=g(x)+\sum_{k\in\mathbb{Z}}\left\{w\int_{\mathbb{R}}\psi(wu-k)\left[\int_x^u\frac{g^{(r)}(t)}{(r-1)!}(u-t)^{r-1}dt\right]\right\}\varphi(wx-k),\quad x\in\mathbb{R},
\]
in view of (i) and (ii). Hence, we can write as follows
\begin{equation*}
\begin{split}
&\left|\left(\mathcal{D}_w^{\varphi,\psi}g\right)(x)-g(x)\right|
\\
&\le\frac{1}{(r-1)!}\sum_{k\in\mathbb{Z}}\left\{w\int_{\mathbb{R}}\left|\psi(wu-k)\right|\left|\int_x^u\left|g^{(r)}(t)\right|\left|u-t\right|^{r-1}dt\right|du\right\}\left|\varphi(wx-k)\right|
\\
&\le\frac{\|g^{(r)}\|_\infty}{(r-1)!}\sum_{k\in\mathbb{Z}}\left\{w\int_{\mathbb{R}}\left|\psi(wu-k)\right|\left|\int_x^u\left|u-t\right|^{r-1}dt\right|du\right\}\left|\varphi(wx-k)\right|
\\
&\le\frac{\|g^{(r)}\|_\infty}{r!}\sum_{k\in\mathbb{Z}}\left\{w\int_{\mathbb{R}}\left|\psi(wu-k)\right|\left|u-x\right|^{r}du\right\}\left|\varphi(wx-k)\right|
\\
&\le\frac{\|g^{(r)}\|_\infty}{w^rr!}\sum_{\nu=0}^r\binom{r}{\nu}\sum_{k\in\mathbb{Z}}\left|k-wx\right|^{r-\nu}\left\{w\int_{\mathbb{R}}\left|\psi(wu-k)\right|\left|wu-k\right|^{\nu}du\right\}\left|\varphi(wx-k)\right|
\\
&=\frac{\|g^{(r)}\|_\infty}{w^rr!}\sum_{\nu=0}^r\binom{r}{\nu}\sum_{k\in\mathbb{Z}}\left|k-wx\right|^{r-\nu}\left|\varphi(wx-k)\right|\widetilde{M}_\nu(\psi)
\\
&\le\frac{\|g^{(r)}\|_\infty}{w^rr!}\sum_{\nu=0}^r\binom{r}{\nu}M_{r-\nu}(\varphi)\widetilde{M}_\nu(\psi)
\\
&=:\frac{\|g^{(r)}\|_\infty}{w^r}\cdot\lambda<+\infty,
\end{split}
\end{equation*}
since $\varphi$ and $\psi$ satisfy ($\Theta$) with $\theta>r+1$ (see Remark \ref{preliminaryremark}).
Thus, by Theorem \ref{thm2} and Remark \ref{rmk1}, we have
\begin{equation*}
\begin{split}
\|\mathcal{D}_w^{\varphi,\psi}f-f\|_\infty&\le(\mu+1)\|f-g\|_\infty+\|\mathcal{D}_w^{\varphi,\psi}g-g\|_\infty
\\
&\le(\mu+1)\left[\|f-g\|_\infty+\frac{\|g^{(r)}\|_\infty\lambda}{w^r(\mu+1)}\right].
\end{split}
\end{equation*}
Finally, by Theorem \ref{Johnen} and the well-known properties of $\omega_r(f,\delta)$, there exists a constant $C_2>0$ such that
\begin{equation*}
\begin{split}
\|\mathcal{D}_w^{\varphi,\psi}f-f\|_\infty&\le(\mu+1)\cdot \mathcal{K}\left(f,w^{-r}\frac{\lambda}{\mu+1},C(\mathbb{R}), C^{r}(\mathbb{R})\right)
\\
&\le(\mu+1)C_2\cdot\omega_r\left(f,w^{-1}\left(\frac{\lambda}{\mu+1}\right)^{\frac{1}{r}}\right)_\infty
\\
&\le(\mu+1)C_2\left[\left(\frac{\lambda}{\mu+1}\right)^{\frac{1}{r}}+1\right]^r\cdot\omega_r\left(f,w^{-1}\right)_\infty
\\
&=: \Lambda_\infty \cdot\omega_r\left(f,\frac{1}{w}\right)_\infty.
\end{split}
\end{equation*}
\end{proof}
As an immediate consequence of both Theorem \ref{Linf} and the well-known inequality
\begin{equation}\label{modulusinequality}
\omega_r(f,\delta)_p\le2^{r-m}\omega_m(f,\delta)_p,\quad\delta>0,
\end{equation}
for $1\le p\le+\infty$ and $1\le m<r$, we can establish the following.
\begin{cor} Under the assumptions of Theorem \ref{Linf}, let $f\in Lip^*(\alpha,L^\infty)$, with $0<\alpha<r$. Thus, we have
\[
\|\mathcal{D}_w^{\varphi,\psi}f-f\|_\infty=\mathcal{O}(w^{-\alpha}),\quad w\to+\infty.
\]
Moreover, if $f\in Lip(r,L^\infty)=C^{r}(\mathbb{R})$, we have
\[
\|\mathcal{D}_w^{\varphi,\psi}f-f\|_\infty=\mathcal{O}(w^{-r}),\quad w\to+\infty.
\]
\end{cor}
\begin{teorema}\label{Lp} Let $1\le p<+\infty$ and $\varphi$ and $\psi$ be the discrete and continuous kernel respectively satisfying $(i)$ and $(ii)$, such that $\varphi$ fulfils condition $(\Theta)$ with $\theta>pr+1$ and $\psi$ is a compactly supported function. Then, for any $f\in L^p(\mathbb{R})$, there exists a constant $\Lambda_p>0$ such that
\[
\|\mathcal{D}_w^{\varphi,\psi}f-f\|_p\le \Lambda_p\cdot\omega_r\left(f,\frac{1}{w}\right)_p,\qquad w>0.
\]
\end{teorema}
\begin{proof}
Let $f\in L^p(\mathbb{R})$ be fixed and $g\in W^{r,p}(\mathbb{R})$. Using (\ref{taylor}) for $g$, and applying Jensen inequality three times, we can obtain what follows
\begin{align}\label{num}
&\|\mathcal{D}_w^{\varphi,\psi}g-g\|_p^p\nonumber
\\
&\le\frac{1}{\left[(r-1)!\right]^p}\int_{\mathbb{R}}\left|\sum_{k\in\mathbb{Z}}\left\{w\int_{\mathbb{R}}\left|\psi(wu-k)\right|\left|\int_x^u\left|g^{(r)}(t)\right|\left|u-t\right|^{r-1}dt\right|du\right\}\left|\varphi(wx-k)\right|\right|^pdx\nonumber
\\
&\le\frac{M_0(\varphi)^{p-1}}{\left[(r-1)!\right]^p}\int_{\mathbb{R}}\sum_{k\in\mathbb{Z}}\left|w\int_{\mathbb{R}}\left|\psi(wu-k)\right|\left|\int_x^u\left|g^{(r)}(t)\right|\left|u-t\right|^{r-1}dt\right|du\right|^p\left|\varphi(wx-k)\right|dx\nonumber
\\
&\le\frac{M_0(\varphi)^{p-1}\|\psi\|_1^{p-1}}{\left[(r-1)!\right]^p}\int_{\mathbb{R}}\sum_{k\in\mathbb{Z}}\left\{w\int_{\mathbb{R}}\left|\psi(wu-k)\right|\left|\int_x^u\left|g^{(r)}(t)\right|\left|u-t\right|^{r-1}dt\right|^pdu\right\}\left|\varphi(wx-k)\right|dx \nonumber
\\
&\le\frac{M_0(\varphi)^{p-1}\|\psi\|_1^{p-1}}{\left[(r-1)!\right]^p}\int_{\mathbb{R}}\sum_{k\in\mathbb{Z}}\left\{w\int_{\mathbb{R}}\left|\psi(wu-k)\right|\left|u-x\right|^{p(r-1)}\left|\int_x^u\left|g^{(r)}(t)\right|dt\right|^pdu\right\}
\\
&\qquad\times\left|\varphi(wx-k)\right|dx\nonumber
\end{align}
\begin{align}
&\le\frac{M_0(\varphi)^{p-1}\|\psi\|_1^{p-1}}{\left[(r-1)!\right]^p}\int_{\mathbb{R}}\sum_{k\in\mathbb{Z}}\left\{w\int_{\mathbb{R}}\left|\psi(wu-k)\right|\left|u-x\right|^{pr-1}\left|\int_x^u\left|g^{(r)}(t)\right|^pdt\right|du\right\}\nonumber
\\
&\qquad\times\left|\varphi(wx-k)\right|dx.\nonumber
\end{align}
We suppose that supp$(\psi)\subseteq[-\sigma,\sigma]$, $\sigma>0$. Thus, we can write what follows
\begin{equation*}
\begin{split}
\|\mathcal{D}_w^{\varphi,\psi}g-g\|_p^p&\le\frac{M_0(\varphi)^{p-1}\|\psi\|_1^{p-1}}{\left[(r-1)!\right]^p}\int_{\mathbb{R}}\sum_{k\in\mathbb{Z}}\left\{w\int_{-\frac{\sigma}{w}+\frac{k}{w}}^{\frac{\sigma}{w}+\frac{k}{w}}\left|\psi(wu-k)\right|\left|u-x\right|^{pr-1}\right.
\\
&\left.\qquad\times\left|\int_x^u\left|g^{(r)}(t)\right|^pdt\right|du\right\}\left|\varphi(wx-k)\right|dx.
\end{split}
\end{equation*}
Now, we can divide the right side of the above inequality in two terms 
\begin{equation*}
\|\mathcal{D}_w^{\varphi,\psi}g-g\|_p^p\le\frac{M_0(\varphi)^{p-1}\|\psi\|_1^{p-1}}{\left[(r-1)!\right]^p}\left\{\mathcal{I}_1+\mathcal{I}_2\right\},
\end{equation*}
where
\begin{equation*}
\begin{split}
&\mathcal{I}_1:=\int_{\mathbb{R}}\sum_{k\in\mathbb{Z}}\left\{w\int_{-\frac{\sigma}{w}+\frac{k}{w}}^{\frac{\sigma}{w}+\frac{k}{w}}\left|\psi(wu-k)\right|\left|u-x\right|^{pr-1}\left|\int_x^u\left|g^{(r)}(t)\right|^pdt-\int_{u+x-\frac{k}{w}}^u\left|g^{(r)}(t)\right|^pdt\right|du\right\}
\\
&\qquad\times\left|\varphi(wx-k)\right|dx
\end{split}
\end{equation*}
and
\begin{equation*}
\begin{split}
&\mathcal{I}_2:=\int_{\mathbb{R}}\sum_{k\in\mathbb{Z}}\left\{w\int_{-\frac{\sigma}{w}+\frac{k}{w}}^{\frac{\sigma}{w}+\frac{k}{w}}\left|\psi(wu-k)\right|\left|u-x\right|^{pr-1}\left|\int_{u+x-\frac{k}{w}}^u\left|g^{(r)}(t)\right|^pdt\right|du\right\}\left|\varphi(wx-k)\right|dx.
\end{split}
\end{equation*}
We first focus on $\mathcal{I}_1$, that can be obviously simplified as follows
\begin{equation*}
\mathcal{I}_1=\int_{\mathbb{R}}\sum_{k\in\mathbb{Z}}\left\{w\int_{-\frac{\sigma}{w}+\frac{k}{w}}^{\frac{\sigma}{w}+\frac{k}{w}}\left|\psi(wu-k)\right|\left|u-x\right|^{pr-1}\left|\int_{x}^{u+x-\frac{k}{w}}\left|g^{(r)}(t)\right|^pdt\right|du\right\}\left|\varphi(wx-k)\right|dx.
\end{equation*}
Putting $y=t-x$ and using Fubini-Tonelli theorem, we get
\begin{equation*}
\begin{split}
\mathcal{I}_1&=\int_{\mathbb{R}}\sum_{k\in\mathbb{Z}}\left\{w\int_{-\frac{\sigma}{w}+\frac{k}{w}}^{\frac{\sigma}{w}+\frac{k}{w}}\left|\psi(wu-k)\right|\left|u-x\right|^{pr-1}\left|\int_{0}^{u-\frac{k}{w}}\left|g^{(r)}(y+x)\right|^pdy\right|du\right\}\left|\varphi(wx-k)\right|dx
\\
&\le\int_{-\frac{\sigma}{w}}^{\frac{\sigma}{w}}\left[\int_{\mathbb{R}}\left|g^{(r)}(y+x)\right|^p\sum_{k\in\mathbb{Z}}\left\{w\int_{-\frac{\sigma}{w}+\frac{k}{w}}^{\frac{\sigma}{w}+\frac{k}{w}}\left|\psi(wu-k)\right|\left|u-x\right|^{pr-1}du\right\}\left|\varphi(wx-k)\right|dx\right]dy.
\end{split}
\end{equation*}
Now, one can observe that
 \begin{equation*}
\begin{split}
&\sum_{k\in\mathbb{Z}}\left\{w\int_{-\frac{\sigma}{w}+\frac{k}{w}}^{\frac{\sigma}{w}+\frac{k}{w}}\left|\psi(wu-k)\right|\left|u-x\right|^{pr-1}du\right\}\left|\varphi(wx-k)\right|
\\
&=\sum_{k\in\mathbb{Z}}\left\{w\sum_{\nu=0}^{pr-1}\binom{pr-1}{\nu}\int_{-\frac{\sigma}{w}+\frac{k}{w}}^{\frac{\sigma}{w}+\frac{k}{w}}\left|\psi(wu-k)\right|\left|u-\frac{k}{w}\right|^{\nu}du\right\}\left|\varphi(wx-k)\right|\left|\frac{k}{w}-x\right|^{pr-1-\nu}
\\
&\le\frac{1}{w^{pr-1}}\sum_{\nu=0}^{pr-1}\binom{pr-1}{\nu}M_{pr-1-\nu}(\varphi)\widetilde{M}_{\nu}(\psi)<+\infty,
\end{split}
\end{equation*}
since $\varphi$ and $\psi$ satisfy ($\Theta$) with $\theta>pr+1$. Hence, we get
\begin{equation*}
\begin{split}
\mathcal{I}_1&\le \frac{1}{w^{pr-1}}\int_{-\frac{\sigma}{w}}^{\frac{\sigma}{w}}\left[\int_{\mathbb{R}}\left|g^{(r)}(y+x)\right|^p\sum_{\nu=0}^{pr-1}\binom{pr-1}{\nu}M_{pr-1-\nu}(\varphi)\widetilde{M}_{\nu}(\psi)dx\right]dy
\\
&=\frac{1}{w^{pr-1}}\left[\sum_{\nu=0}^{pr-1}\binom{pr-1}{\nu}M_{pr-1-\nu}(\varphi)\widetilde{M}_{\nu}(\psi)\right]\int_{-\frac{\sigma}{w}}^{\frac{\sigma}{w}}\|g^{(r)}(y+\cdot)\|_p^pdy
\\
&=\frac{1}{w^{pr}}\left[\sum_{\nu=0}^{pr-1}\binom{pr-1}{\nu}M_{pr-1-\nu}(\varphi)\widetilde{M}_{\nu}(\psi)\right]2\sigma\|g^{(r)}\|_p^p
\\
&=:\frac{1}{w^{pr}}\|g^{(r)}\|_p^p\cdot\lambda_1<+\infty.
\end{split}
\end{equation*}

We now consider $\mathcal{I}_2$. By the change of variable $y=t-u$, the convexity of $|\cdot|^{pr-1}$ and the application of the Fubini-Tonelli theorem, we have
\begin{equation*}
\begin{split}
\mathcal{I}_2&=\int_{\mathbb{R}}\sum_{k\in\mathbb{Z}}\left\{w\int_{-\frac{\sigma}{w}+\frac{k}{w}}^{\frac{\sigma}{w}+\frac{k}{w}}\left|\psi(wu-k)\right|\left|u-x\right|^{pr-1}\left|\int_{x-\frac{k}{w}}^0\left|g^{(r)}(y+u)\right|^pdy\right|du\right\}\left|\varphi(wx-k)\right|dx
\\
&\le 2^{pr-2}\int_{\mathbb{R}}\sum_{k\in\mathbb{Z}}\left\{w\int_{-\frac{\sigma}{w}+\frac{k}{w}}^{\frac{\sigma}{w}+\frac{k}{w}}\left|\psi(wu-k)\right|\left[\left(\frac{\sigma}{w}\right)^{pr-1}+\left|\frac{k}{w}-x\right|^{pr-1}\right]\right.
\\
&\left.\qquad\times\left|\int_{x-\frac{k}{w}}^0\left|g^{(r)}(y+u)\right|^pdy\right|du\right\}\left|\varphi(wx-k)\right|dx
\\
&\le 2^{pr-2}\sum_{k\in\mathbb{Z}}\int_{\mathbb{R}}\left\{w\int_{-\frac{\sigma}{w}+\frac{k}{w}}^{\frac{\sigma}{w}+\frac{k}{w}}\left|\psi(wu-k)\right|\left[\left(\frac{\sigma}{w}\right)^{pr-1}+\left|\frac{k}{w}-x\right|^{pr-1}\right]\right.
\\
&\left.\qquad\times\left|\int_{x-\frac{k}{w}}^0\left|g^{(r)}(y+u)\right|^pdy\right|du\right\}\left|\varphi(wx-k)\right|dx.
\end{split}
\end{equation*}
Now, setting $z=wx-k$ and applying again the Fubini-Tonelli theorem, we obtain
\begin{equation*}
\begin{split}
\mathcal{I}_2&\le 2^{pr-2}\sum_{k\in\mathbb{Z}}\int_{\mathbb{R}}\left\{\int_{-\frac{\sigma}{w}+\frac{k}{w}}^{\frac{\sigma}{w}+\frac{k}{w}}\left|\psi(wu-k)\right|\left[\left(\frac{\sigma}{w}\right)^{pr-1}+\left|\frac{z}{w}\right|^{pr-1}\right]\right.
\\
&\left.\qquad\times\left|\int_{x-\frac{k}{w}}^0\left|g^{(r)}(y+u)\right|^pdy\right|du\right\}\left|\varphi(z)\right|dz
\\
&\le 2^{pr-2}\int_{\mathbb{R}}\left|\varphi(z)\right|\sum_{k\in\mathbb{Z}}\left\{\int_{-\frac{\sigma}{w}+\frac{k}{w}}^{\frac{\sigma}{w}+\frac{k}{w}}\left|\psi(wu-k)\right|\left[\left(\frac{\sigma}{w}\right)^{pr-1}+\left|\frac{z}{w}\right|^{pr-1}\right]\right.
\\
&\left.\qquad\times\left(\int_{|y|\le|z|/w}\left|g^{(r)}(y+u)\right|^pdy\right)du\right\}dz
\\
&=\frac{ 2^{pr-2}}{w^{pr-1}}\int_{\mathbb{R}}\left(\sigma^{pr-1}+|z|^{pr-1}\right)\left|\varphi(z)\right|
\\
&\qquad\times\left(\int_{|y|\le|z|/w}\sum_{k\in\mathbb{Z}}\left\{\int_{-\frac{\sigma}{w}+\frac{k}{w}}^{\frac{\sigma}{w}+\frac{k}{w}}\left|\psi(wu-k)\right|\left|g^{(r)}(y+u)\right|^pdu\right\}dy\right)dz
\end{split}
\end{equation*}
\begin{equation*}
\begin{split}
&\le\frac{ 2^{pr-2}M_0(\psi)}{w^{pr-1}}\int_{\mathbb{R}}\left(\sigma^{pr-1}+|z|^{pr-1}\right)\left|\varphi(z)\right|\left(\int_{|y|\le|z|/w}\|g^{(r)}(y+\cdot)\|^p_p dy\right)dz
\\
&\le\frac{ 2^{pr-1}M_0(\psi)}{w^{pr}}\|g^{(r)}\|_p^p\int_{\mathbb{R}}\left(\sigma^{pr-1}|z|+|z|^{pr}\right)\left|\varphi(z)\right|dz
\\
&\le\frac{ 2^{pr-1}M_0(\psi)}{w^{pr}}\|g^{(r)}\|_p^p\left\{\sigma^{pr-1}\widetilde{M}_1(\varphi)+\widetilde{M}_{pr}(\varphi)\right\}<+\infty,
\end{split}
\end{equation*}
since $\varphi$ and $\psi$ satisfy ($\Theta$) with $\theta>pr+1$.\\
Now, we are able to summarize all the previous estimates as follows

\begin{equation*}
\begin{split}
&\|\mathcal{D}_w^{\varphi,\psi}g-g\|_p
\\
&\le w^{-r}\|g^{(r)}\|_p\cdot \frac{M_0(\varphi)^{\frac{p-1}{p}}\|\psi\|_1^{\frac{p-1}{p}}}{(r-1)!}\left(\lambda_1+2^{pr-1}M_0(\psi)\left\{\sigma^{pr-1}\widetilde{M}_1(\varphi)+\widetilde{M}_{pr}(\varphi)\right\}\right)^{\frac{1}{p}}
\\
&=:w^{-r}\|g^{(r)}\|_p\cdot\lambda_2.
\end{split}
\end{equation*}
Therefore, by Theorem \ref{thm2} we have
\begin{equation*}
\begin{split}
&\|\mathcal{D}_w^{\varphi,\psi}f-f\|_p\le(\mu+1)\left[\|f-g\|_p+\frac{\|g^{(r)}\|_p\lambda_2}{w^r(\mu+1)}\right].
\end{split}
\end{equation*}
Finally, passing to the infimum with respect to $g\in W^{r,p}(\mathbb{R})$ and using Theorem \ref{Johnen}, we have
\begin{equation*}
\begin{split}
\|\mathcal{D}_w^{\varphi,\psi}f-f\|_p&\le(\mu+1)\cdot\mathcal{K}\left(f,w^{-r}\frac{\lambda_2}{\mu+1},L^p,W^{r,p}\right)
\\
&=: \Lambda_p\cdot\omega_r\left(f,w^{-1}\right)_p,
\end{split}
\end{equation*}
where $\Lambda_p>0$ is a suitable constant obtained proceeding as at the end of the proof of Theorem \ref{Linf}.
\end{proof}

For the sake of generality, we now want to avoid the assumption for $\psi$ to be a compactly supported function. To this aim, we use the well-known \textit{Hardy–Littlewood maximal function} \cite{ACP2012}, defined as

\[
	\mathcal{M}f(x):=\sup_{\stackrel{u\in\mathbb{R}}{u \neq x} }\frac{1}{\left|x-u\right|}\int_{x}^u\left|f(t)\right|dt,\qquad x\in \mathbb{R},
\]
for a locally integrable function $f:\mathbb{R}\to\mathbb{R}$. The celebrated theorem of Hardy, Littlewood and Wiener asserts that $\mathcal{M}f$ is bounded on $L^p(\mathbb{R})$ for $1<p\le+\infty$ for every $f\in L^p(\mathbb{R})$, namely
\begin{equation}\label{HLM}
\|\mathcal{M}f\|_p\le C_p\|f\|_p,
\end{equation}
where the constant $C_p$ depends only on $p$ (Theorem I.1 of \cite{Stein1970}). However, the corresponding result for $p=1$ fails.
\begin{teorema}\label{thmhl} Let $\varphi$ and $\psi$ be the discrete and the continuous kernel respectively satisfying $(i)$, $(ii)$, and condition $(\Theta)$ with $\theta>pr+1$. Then, for any $f\in L^p(\mathbb{R})$, with $1<p<+\infty$, there exists a constant $\Gamma_p>0$ such that
\[
\|\mathcal{D}_w^{\varphi,\psi}f-f\|_p\le \Gamma_p\cdot\omega_r\left(f,\frac{1}{w}\right)_p,\qquad w>0.
\]
\end{teorema}
\begin{proof}
Let $f\in L^p(\mathbb{R})$ be fixed and $g\in W^{r,p}(\mathbb{R})$. Using (\ref{num}), we have
\begin{equation*}
\begin{split}
&\|\mathcal{D}_w^{\varphi,\psi}g-g\|_p^p
\\
&\le\frac{M_0(\varphi)^{p-1}\|\psi\|_1^{p-1}}{\left[(r-1)!\right]^p}\int_{\mathbb{R}}\sum_{k\in\mathbb{Z}}\left\{w\int_{\mathbb{R}}\left|\psi(wu-k)\right|\left|u-x\right|^{pr}\left|\int_x^u\frac{\left|g^{(r)}(t)\right|}{\left|u-x\right|}dt\right|^pdu\right\}\left|\varphi(wx-k)\right|dx
\\
&\le\frac{M_0(\varphi)^{p-1}\|\psi\|_1^{p-1}}{\left[(r-1)!\right]^p}\int_{\mathbb{R}}\sum_{k\in\mathbb{Z}}\left\{w\int_{\mathbb{R}}\left|\psi(wu-k)\right|\left|u-x\right|^{pr}\left|\mathcal{M}g^{(r)}(x)\right|^pdu\right\}\left|\varphi(wx-k)\right|dx.
\end{split}
\end{equation*}
Thus, we can estimate the right-hand side as follows
\begin{equation*}
\begin{split}
&\|\mathcal{D}_w^{\varphi,\psi}g-g\|_p^p
\\
&\le\frac{M_0(\varphi)^{p-1}\|\psi\|_1^{p-1}}{\left[(r-1)!\right]^p}\int_{\mathbb{R}}\left|\mathcal{M}g^{(r)}(x)\right|^p\sum_{k\in\mathbb{Z}}\left\{w\int_{\mathbb{R}}\left|\psi(wu-k)\right|\left|u-x\right|^{pr}du\right\}\left|\varphi(wx-k)\right|dx
\\
&=\frac{M_0(\varphi)^{p-1}\|\psi\|_1^{p-1}}{w^{pr}\left[(r-1)!\right]^p}\int_{\mathbb{R}}\left|\mathcal{M}g^{(r)}(x)\right|^p
\\
&\qquad\times\sum_{k\in\mathbb{Z}}\left\{w\sum_{\nu=0}^{pr}\binom{pr}{\nu}\int_{\mathbb{R}}\left|\psi(wu-k)\right|\left|wu-k\right|^{\nu}du\right\}\left|k-wx\right|^{pr-\nu}\left|\varphi(wx-k)\right|dx
\\
&\le\frac{M_0(\varphi)^{p-1}\|\psi\|_1^{p-1}}{w^{pr}\left[(r-1)!\right]^p}\int_{\mathbb{R}}\left|\mathcal{M}g^{(r)}(x)\right|^p\left\{\sum_{\nu=0}^{pr}\binom{pr}{\nu}\widetilde{M}_\nu(\psi)M_{pr-\nu}(\varphi)\right\}dx
\\
&=:\frac{\lambda}{w^{pr}}\cdot\|\mathcal{M}g^{(r)}\|_p^p<+\infty, 
\end{split}
\end{equation*}
being $g^{(r)}\in L^p(\mathbb{R})$ and $\lambda<+\infty$ since condition ($\Theta$) is satisfied for $\theta>pr+1$. Now, by (\ref{HLM}) there exists a constant $C_p>0$ such that
\begin{equation*}
\begin{split}
&\|\mathcal{D}_w^{\varphi,\psi}g-g\|_p^p
\le\frac{\lambda}{w^{pr}} C_p\|g^{(r)}\|_p^p<+\infty.
\end{split}
\end{equation*}
Hence
\begin{equation*}
\begin{split}
\|\mathcal{D}_w^{\varphi,\psi}g-g\|_p&\le w^{-r}\|g^{(r)}\|_p\left(\lambda C_p\right)^{\frac{1}{p}}
=:w^{-r}\|g^{(r)}\|_p\cdot\lambda_1.
\end{split}
\end{equation*}
Finally, we get the thesis as in the previous cases, obtaining
\begin{equation*}
\begin{split}
\|\mathcal{D}_w^{\varphi,\psi}f-f\|_p
\le \Gamma_p\cdot\omega_r\left(f,w^{-1}\right)_p,
\end{split}
\end{equation*}
for a suitable constant $\Gamma_p>0$.
\end{proof}

Now, we may provide the following corollary, by using again inequality (\ref{modulusinequality}) together with Theorem \ref{Lp} and Theorem \ref{thmhl}.

\begin{cor} Let $1\leq p<+\infty$ and $f\in Lip^*(\alpha,L^p)$, with $0<\alpha<r$. 
\begin{description}
\item[{(a)}] If $1<p<+\infty$, under the assumptions of Theorem \ref{thmhl}, we have
\[
\|\mathcal{D}_w^{\varphi,\psi}f-f\|_p=\mathcal{O}(w^{-\alpha}),\quad w\to+\infty.
\]
In particular, if $f\in Lip(r,L^p)=W^{r,p}(\mathbb{R})$, there holds
\[
\|\mathcal{D}_w^{\varphi,\psi}f-f\|_p=\mathcal{O}(w^{-r}),\quad w\to+\infty.
\]
\item[{(b)}] If $p=1$, under the assumptions of Theorem \ref{Lp}, we have
\[
\|\mathcal{D}_w^{\varphi,\psi}f-f\|_1=\mathcal{O}(w^{-\alpha}),\quad w\to+\infty.
\]
Moreover, if $f\in Lip(r,L^1)= W^{r-1}(BV)\cap L^1(\mathbb{R})$, there holds
\[
\|\mathcal{D}_w^{\varphi,\psi}f-f\|_1=\mathcal{O}(w^{-r}),\quad w\to+\infty.
\]
\end{description}
\end{cor}

\section[Inverse theorem of approximation and characterization of generalized Lipschitz classes]{An inverse theorem of approximation and a characterization of the generalized Lipschitz classes}
In order to achieve an inverse approximation result, we premise the following useful Lemma (see Lemma 5.3 of \cite{2022JFAA}).
\begin{lemma}[]\label{moment} Let $\varphi\in W^{r,1}(\mathbb{R})$ be a discrete kernel satisfying $(i)$ and $(\Theta)$ with $\theta>r+1$. Then
\[
m_\mu(\varphi^{(\nu)},u)=\begin{cases}
			0, & \text{if $\nu\ne \mu$}\\
            \nu!, & \text{if $\nu=\mu$,}
		 \end{cases}
\qquad u\in\mathbb{R},\]
where $\nu=1,\dots,r$ and $\mu=0,\dots,\nu$.
\end{lemma}
Now, we are able to prove the following inversion theorem.
\begin{teorema}\label{invthm} Let $\varphi\in W^{r,1}(\mathbb{R})$ be a discrete kernel satisfying $(i)$ and $\psi\in L^1(\mathbb{R})$ be a continuous kernel.
Let $f\in L^p(\mathbb{R})$, $1\le p\le+\infty$ such that
\begin{equation}\label{convergence}
\|\mathcal{D}_w^{\varphi,\psi}f-f\|_p=\mathcal{O}(w^{-\alpha}),
\end{equation}
as $w\to+\infty$, with $0<\alpha<r$. Then
\begin{description}
\item[{(a)}] if $p=+\infty$, and $\varphi^{(i)}$, $i=1,\dots, r$, and $\psi$ satisfy $(\Theta)$ for $\theta>r+1$, it turns out that
\[
f\in Lip^*(\alpha,L^\infty);
\]
\item[{(b)}] if $1\le p<+\infty$, and $\varphi^{(i)}$, $i=1,\dots, r$, and $\psi$ satisfy $(\Theta)$ for $\theta>pr+1$ (if $p=1$, $\psi$ is assumed to be a compactly supported function), it turns out that
\[
f\in Lip^*(\alpha,L^p).
\]
\end{description}
\end{teorema}
\begin{proof}
Let $1\le p\le+\infty$. By (\ref{convergence}), there exist $C,\overline{w}>0$ such that
\[
\|\mathcal{D}_w^{\varphi,\psi}f-f\|_p\le \dfrac{C}{w^\alpha},
\]
for every $w\ge\overline{w}$. We set $\overline{\delta}:=\dfrac{1}{\overline{w}}$ and we now consider $0<\delta<\overline{\delta}$, $w\ge\overline{w}$. Being  $0<\alpha<r$, we can fix $m\in\mathbb{N}$ such that $1\le m\le r$ so that $m=1$ if $0<\alpha<1$ or $2\le m\le r$ if $m-1\le\alpha<m$. In Theorem 3.3 of \cite{DurrRegularization}, we prove that $\mathcal{D}_w^{\varphi,\psi}f\in W^{r,p}(\mathbb{R})$. Let now $g\in W^{r,p}(\mathbb{R})$. From (\ref{modulusinequality}), it turns out that
\begin{equation*}
\begin{split}
\omega_m(f,\delta)_p&\le\omega_m(f-\mathcal{D}_w^{\varphi,\psi}f,\delta)_p+\omega_m(\mathcal{D}_w^{\varphi,\psi}f,\delta)_p
\\
&\le 2^m\|\mathcal{D}_w^{\varphi,\psi}f-f\|_p+\dfrac{1}{C_1}\mathcal{K}\left(\mathcal{D}_w^{\varphi,\psi}f,\delta^m,L^p,W^{r,p}\right)
\\
&\le \dfrac{2^m C}{w^\alpha}+\dfrac{\delta^m}{C_1}\|\left(\mathcal{D}_w^{\varphi,\psi}f\right)^{(m)}\|_p
\\
&\le \dfrac{2^m C}{w^\alpha}+\dfrac{\delta^m}{C_1}\left\{\lVert\left(\mathcal{D}_w^{\varphi,\psi}f\right)^{(m)}-\left(\mathcal{D}_w^{\varphi,\psi}g\right)^{(m)}\rVert_p+\|\left(\mathcal{D}_w^{\varphi,\psi}g\right)^{(m)}\|_p\right\}
\\
&\le \dfrac{2^m C}{w^\alpha}+\dfrac{\delta^m}{C_1}\left\{w^m K\|f-g\|_p+\|\left(\mathcal{D}_w^{\varphi,\psi}f\right)^{(m)}\|_p\right\},
\end{split}
\end{equation*}
where
\[
K=\begin{cases} M_0(\varphi^{(m)})^{\frac{p-1}{p}}\|\psi\|_1^{\frac{p-1}{p}}\|\varphi^{(m)}\|_1^{\frac{1}{p}}M_0(\psi)^{\frac{1}{p}}, & \mbox{if }\mbox{ $1\le p<+\infty$} \\ M_0(\varphi^{(m)})\|\psi\|_1, & \mbox{if }\mbox{ $p=+\infty$,}
\end{cases}
\]
arguing as in Theorem \ref{thm2}.
\\
By Lemma \ref{moment}, we know that $m_0(\varphi^{(m)},u)=0$, $u\in\mathbb{R}$, and then by Theorem 3.3 of \cite{DurrRegularization}, we have
\begin{equation*}
\begin{split}
\left(\mathcal{D}_w^{\varphi,\psi}g\right)^{(m)}(x)&=w^m\sum_{k\in\mathbb{Z}}\left\{w\int_{\mathbb{R}}\psi(wu-k)g(u)du\right\}\varphi^{(m)}(wx-k)
\\
&=w^m\sum_{k\in\mathbb{Z}}\left\{w\int_{\mathbb{R}}\psi(wu-k)g(u)du\right\}\varphi^{(m)}(wx-k)-g(x)w^m\sum_{k\in\mathbb{Z}}\varphi^{(m)}(wx-k)
\\
&=w^m\sum_{k\in\mathbb{Z}}\left\{w\int_{\mathbb{R}}\psi(wu-k)\left[g(u)-g(x)\right]du\right\}\varphi^{(m)}(wx-k),\quad x\in\mathbb{R}.
\end{split}
\end{equation*}
We now expand the operator $\mathcal{D}_w^{\varphi,\psi}g^{(m)}$ as in (\ref{taylor}), obtaining
\begin{equation*}
\begin{split}
\left(\mathcal{D}_w^{\varphi,\psi}g\right)^{(m)}(x)&=w^m\sum_{i=1}^{m-1}\dfrac{g^{(i)}(x)}{i!}\sum_{\nu=0}^i\binom{i}{\nu}\sum_{k\in\mathbb{Z}}\left(\dfrac{k}{w}-x\right)^{i-\nu}
\\
&\qquad\times\left\{w\int_{\mathbb{R}}\psi(wu-k)\left(u-\dfrac{k}{m}\right)^\nu\right\}\varphi^{(m)}(wx-k)
\\
&\qquad+w^m\sum_{k\in\mathbb{Z}}\left\{w\int_\mathbb{R}\psi(wu-k)\left[\int_x^u\frac{g^{(m)}(t)}{(m-1)!}(u-t)^{m-1}dt\right]du\right\}\varphi^{(m)}(wx-k)
\end{split}
\end{equation*}
\begin{equation*}
\begin{split}
&=w^m\sum_{i=1}^{m-1}\dfrac{g^{(i)}(x)}{i!w^i}\sum_{\nu=0}^i\binom{i}{\nu}m_{i-\nu}(\varphi^{(m)},wx)\widetilde{m}_\nu(\psi)
\\
&\qquad+w^m\sum_{k\in\mathbb{Z}}\left\{w\int_\mathbb{R}\psi(wu-k)\left[\int_x^u\frac{g^{(m)}(t)}{(m-1)!}(u-t)^{m-1}dt\right]du\right\}\varphi^{(m)}(wx-k),
\end{split}
\end{equation*}
$x\in\mathbb{R}$. By using again Lemma \ref{moment}, we have that
\[
m_{i-\nu}(\varphi^{(m)},wx)=0,\quad x\in\mathbb{R},
\]
for every $i=1,\dots,m-1$ and $\nu=0,\dots,i$. Hence, we reduce to consider
\begin{equation*}
\begin{split}
\left(\mathcal{D}_w^{\varphi,\psi}g\right)^{(m)}(x)=w^m\sum_{k\in\mathbb{Z}}\left\{w\int_\mathbb{R}\psi(wu-k)\left[\int_x^u\frac{g^{(m)}(t)}{(m-1)!}(u-t)^{m-1}dt\right]du\right\}\varphi^{(m)}(wx-k).
\end{split}
\end{equation*}
Since $\varphi^{(m)}$ and $\psi$ satisfy condition ($\Theta$) for a suitable $\theta>0$, we can approach the same method used in Theorem \ref{Linf} when $p=+\infty$, in Theorem \ref{thmhl} when $1< p<+\infty$, or in Theorem \ref{Lp} when $p=1$, to obtain
\begin{equation*}
\begin{split}
\omega_m(f,\delta)_p&\le \dfrac{2^m C}{w^\alpha}+\dfrac{\delta^mw^m}{C_1}\left\{K\|f-g\|_p+\tau w^{-m}\|g^{(m)}\|_p\right\}
\\
&\le\dfrac{2^m C}{w^\alpha}+\dfrac{\delta^mw^mK}{C_1}\left\{\|f-g\|_p+\dfrac{\tau w^{-m}}{K}\|g^{(m)}\|_p\right\},
\end{split}
\end{equation*}
where $\tau>0$ is a suitable constant.
Passing to the infimum with respect to $g\in W^{r,p}(\mathbb{R})$, we have
\begin{equation}\label{eqinv}
\begin{split}
\omega_m(f,\delta)_p&\le\dfrac{2^m C}{w^\alpha}+\dfrac{\delta^mw^m K}{C_1}\cdot\mathcal{K}\left(f,\dfrac{\tau w^{-m}}{K},L^p,W^{r,p}\right)
\\
&\le\dfrac{K_1}{w^\alpha}+K_2\delta^mw^m\cdot\omega\left(f,\dfrac{1}{w}\right)_p,
\end{split}
\end{equation}
where $K_1$, $K_2$ are two suitable positive constants.
\\
We now fix $A>\bar{w}>1$ sufficiently large so that $2K_2<A^{m-\alpha}$, with $m-\alpha>0$. We take
\[
M:=\max\left\{\omega_m(f,1)_p,2K_1A^\alpha\right\}
\]
and now we are going to show that 
\begin{equation}\label{induction}
\omega_m(f,\delta_n)_p\le M\delta_n^\alpha,
\end{equation}
where $\delta_n:=A^{-n}$, $n\in\mathbb{N}$. If $n=0$, inequality (\ref{induction}) is trivially satisfied, because $\omega_m(f,1)_p\le M\delta_0^\alpha=M$.
For $n>0$, we argue by induction, supposing that (\ref{induction}) holds for $n-1\ge 1$. In (\ref{eqinv}), we take $w=A^{n-1}$ and $\delta=\delta_n$. Now, by the general consideration that, if $a\le b+c$, then either $a\le2b$ or $a\le 2c$, by (\ref{eqinv}) we can get
\begin{equation*}
\begin{split}
\omega_m(f,\delta_n)_p&\le2K_1\delta^\alpha_{n-1}
=2K_1A^\alpha\delta_n^\alpha
\le M\delta_n^\alpha
\end{split}
\end{equation*}
or
\begin{equation*}
\begin{split}
\omega_m(f,\delta_n)_p&\le2K_2A^{m(n-1)}\delta_n^m\omega_m(f,\delta_{n-1})_p
\\
&\le2K_2A^{-m}M\delta_{n-1}^\alpha
\\
&=2K_2A^{\alpha-m}M\delta_n^\alpha
\\
&
<M\delta^\alpha_n,
\end{split}
\end{equation*}
noting that $A^{-\alpha}\delta_{n-1}^\alpha=\delta_n^\alpha$.
\\
Finally, in order to show that $\omega_m(f,\delta)_p=\mathcal{O}(\delta^\alpha)$ as $\delta\to 0^+$, we have to find a suitable constant $N>0$ such that $\omega_m(f,\delta)_p\le N\delta^\alpha$, for every $0<\delta<\delta_1$, with $\delta_1\le\bar{\delta}$. Then, we fix $0<\delta<\delta_1$ and $n\in\mathbb{N}$ such that $\delta_n<\delta\le\delta_{n-1}$. Hence
\begin{equation*}
\begin{split}
\omega_m(f,\delta)_p&\le\omega_m(f,\delta_{n-1})_p
\le M\delta_{n-1}^\alpha
\\
&=MA^\alpha\delta_n^\alpha
=:N\delta_n^\alpha
<N\delta^\alpha.
\end{split}
\end{equation*}
This completes the proof.
\end{proof}

Finally, we are prepared to characterize the generalized Lipschitz classes by the rate of convergence of Durrmeyer sampling operators within the $L^p$-setting. To achieve this, we will employ the direct results outlined in Section \ref{directsection}, together with the inverse approximation theorem established above (Theorem \ref{invthm}).
\begin{teorema}\label{carthm} Let $\varphi\in W^{r,1}(\mathbb{R})$ and $\psi\in L^1(\mathbb{R})$ be a discrete and a continuous kernel, respectively, satisfying $(i)$ and $(ii)$. For any $f\in L^p(\mathbb{R})$, with $1\le p\le+\infty$, and $0<\alpha<r$, it turns out that
\begin{description}
\item[{(I)}] if $p=1$ and $\varphi$, $\varphi^{(i)}$ with $i=1,\dots, r$ satisfy $(\Theta)$ for $\theta>r+1$ and $\psi$ is a compactly supported function, then
\[
f\in Lip^*(\alpha,L^1)
\]
if and only if
\[
\|\mathcal{D}_w^{\varphi,\psi}f-f\|_1=\mathcal{O}(w^{-\alpha}),\qquad w\to+\infty;
\]
\item[{(II)}] if $1< p<+\infty$ and $\varphi$, $\varphi^{(i)}$ with $i=1,\dots, r$ and $\psi$ satisfy $(\Theta)$ for $\theta>pr+1$, then
\[
f\in Lip^*(\alpha,L^p)
\]
if and only if
\[
\|\mathcal{D}_w^{\varphi,\psi}f-f\|_p=\mathcal{O}(w^{-\alpha}),\qquad w\to+\infty;
\]
\item[{(III)}] if $p=+\infty$ and $\varphi$, $\varphi^{(i)}$ with $i=1,\dots, r$ and $\psi$ satisfy $(\Theta)$ for $\theta>r+1$, then
\[
f\in Lip^*(\alpha,L^\infty)
\]
if and only if
\[
\|\mathcal{D}_w^{\varphi,\psi}f-f\|_\infty=\mathcal{O}(w^{-\alpha}),\qquad w\to+\infty.
\]
\end{description}
\end{teorema}

\section{Some applications}
In this section, we want to discuss different instances of kernels $\varphi$ and $\psi$ such that the fundamental assumptions (i) and (ii) are satisfied. Therefore, especially to verify the vanishing moment condition, it is essential to establish their algebraic moments, both the discrete as the continuous versions. First of all, we recall that the \textit{Fourier transform} of a function $f\in L^1(\mathbb{R})$ is defined by
\[
\widehat{f}(v)=\int_{\mathbb{R}}f(u)e^{-iuv}du,\qquad v\in\mathbb{R}.
\]

\begin{figure*}[tbph]
\centering
 \includegraphics[width=10cm]{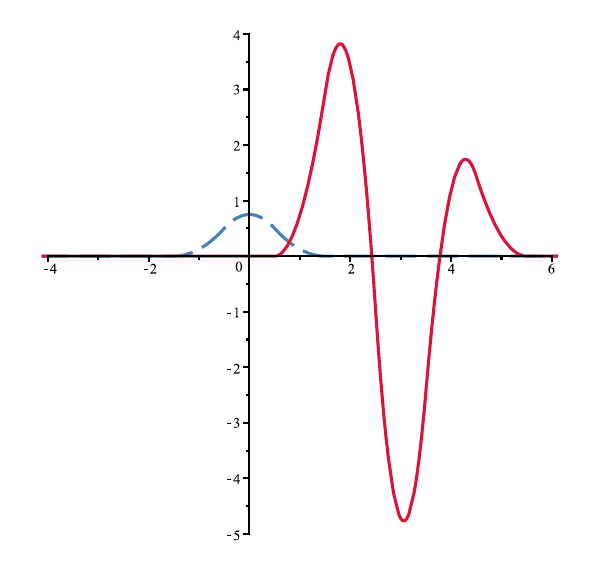}
\caption{The graphs of the kernels $\sigma_3$ (dashed blue line) and $\tau$ (solid red line).}
\label{spline_traslate}
\end{figure*}

A first example can be given considering the well-known \textit{central B-splines of order $n\in\mathbb{N}$}, defined by
\begin{equation}\label{spline}
\sigma_n(u):=\frac{1}{(n-1)!}\sum_{j=0}^n(-1)^j\binom{n}{j}\left(\frac{n}{2}+u-j\right)_+^{n-1},\qquad u\in\mathbb{R},
\end{equation}
where $(\cdot)_+$ denotes the positive part, i.e., $(u)_+:=\max\left\{u,0\right\}$ (see \cite{1981jws}). The corresponding Fourier transform is given by $\displaystyle\widehat{\sigma}_n(v)=\text{sinc}^n\left(\frac{v}{2\pi}\right)$, $v\in\mathbb{R}$, where the sinc-function is defined as $\sin(\pi v)/\pi v$, if $v\ne 0$, and $1$, if $v=0$.
\\
To establish the algebraic discrete moments, namely $m_\nu(\sigma_n, x)$, $x\in\mathbb{R}$, with  $\nu\in\mathbb{N}$, we use the well-known \textit{Poisson summation formula} (PSF) (see \cite{FAA71}). We first observe that for $m\in\mathbb{Z}$ ($m\neq 0$) and $\nu=1,\dots,n-1$, we have $\widehat{\sigma}_n^{(\nu)}(2m\pi)=0$. Therefore, it follows that
\[
m_\nu(\sigma_n,x)=\begin{cases} 0, & \mbox{if }\nu \mbox{ is odd,} \\ (-1)^{\nu/2}\widehat{\sigma}_n^{(\nu)}(0), & \mbox{if } \nu \mbox{ is even},
\end{cases}\qquad \nu=1,\dots,n-1,\ x\in\mathbb{R}.
\]
As to the continuous moments, by differentiating under the sign of integral, we obtain
\begin{equation*}\label{diffsignintegral}
\widehat{\sigma}_n^{(\nu)}(v)=(-i)^\nu\int_{-n/2}^{n/2}\sigma_n(u)u^\nu e^{-ivu}du,\qquad v\in\mathbb{R}.
\end{equation*}
In particular, for $v=0$ we have $\widetilde{m}_\nu(\sigma_n)=(-1)^{\nu/2}\widehat{\sigma}_n^{(\nu)}(0)$, for every $\nu\in\mathbb{N}$. Thus, $\widetilde{m}_\nu(\sigma_n)=m_\nu(\sigma_n, x)$, for every $\nu=1,\dots,n-1$, $x\in\mathbb{R}$. Moreover, a further consequence of PSF is that $m_0(\sigma_n, x)=\widetilde{m}_0(\sigma_n)=1$, for every $n\in\mathbb{N}$, $x\in\mathbb{R}$. We also stress that all the absolute moments, namely $M_\nu(\sigma_n)$ and $\widetilde{M}_\nu(\sigma_n)$, are trivially finite for every $\nu\ge0$ (see also Remark \ref{preliminaryremark}).
\\
In view of the analysis, taking $\varphi=\psi=\sigma_n$, with $n\ge2$, assumptions (i) and (ii) are satisfied for $r=2$, while condition ($\Theta$) is trivially fulfilled for every $\theta>0$. 
Therefore, from Theorem \ref{Linf}, Theorem \ref{Lp} and the corresponding corollaries, we may state the following.
\begin{cor}\label{corex} Let $n\ge 2$. There exists $C>0$ such that
for any $f\in L^p(\mathbb{R})$ with $1\le p\le+\infty$, it turns out that \[\|\mathcal{D}_w^{\sigma_n,\sigma_n}f-f\|_\infty\le C\cdot\omega_2\left(f,\frac{1}{w}\right)_p.\]
In particular, if $f\in Lip(2,L^p)$, we have
\[
\|\mathcal{D}_w^{\sigma_n,\sigma_n}f-f\|_p=\mathcal{O}(w^{-2}),\quad w\to+\infty.
\]
\end{cor}
\noindent We recall that $Lip(2,L^p)=W^{2,p}(\mathbb{R})$ if $p>1$, otherwise $Lip(2,L^1)=W^{1}(BV)\cap L^1(\mathbb{R})$. In particular, if $p=+\infty$, then $W^{2,\infty}(\mathbb{R})=C^2(\mathbb{R})$.
\\
For completeness, we can also characterize the aforementioned generalized Lipschitz classes using Theorem \ref{carthm}. Notably, all the necessary assumptions hold when we consider B-spline kernels $\sigma_n$, with $n\ge 2$.
\begin{cor} Let $n\ge 2$. For any $f\in L^p(\mathbb{R})$, with $1\le p\le +\infty$, it turns out that
\[
f\in Lip^*(\alpha, L^p)
\]
if and only if
\[
\|\mathcal{D}_w^{\sigma_n,\sigma_n}f-f\|_p=\mathcal{O}(w^{-\alpha}),
\]
as $w\to+\infty$, with $0<\alpha<2$.
\end{cor}
We point out that Corollary \ref{corex} holds even if we choose a not necessarily compactly supported kernels, as a consequence of Theorem \ref{thmhl}. In fact, we may take into account the well-known \textit{Jackson kernel of order $N\in\mathbb{N}$}, defined by
\[
J_{N,\alpha}(u):=c_{N,\alpha}\text{sinc}^{2N}\left(\frac{u}{2N\pi\alpha}\right),\qquad u\in\mathbb{R},
\]
with $N\in\mathbb{N}$, $\alpha\ge 1$, and $c_{N,\alpha}$ a non-zero normalization coefficient, given by
\[
c_{N,\alpha}:=\left[\int_\mathbb{R}\text{sinc}^{2N}\left(\frac{u}{2N\pi\alpha}\right)\,du\right]^{-1}.
\]
In particular, it is well-known that $J_{N,\alpha}$ has Fourier transform with support contained in $[-1/\alpha,1/\alpha]$, and by PSF we have $m_0(J_{N,\alpha})=1$ and $m_1(J_{N,\alpha})=-i\widehat{J}_{N,\alpha}^{(1)}(0)=0$. Moreover, being $J_{N,\alpha}$ a positive and even function on $\mathbb{R}$, it turns out that $\widetilde{m}_\nu(J_{N,\alpha})=0$, for every $\nu$ odd. Therefore, the vanishing moment condition (ii) is again satisfied for $r=2$. In addition, it is easy to observe that $J_{N,\alpha}(u)=\mathcal{O}(|u|^{-2N})$, as $|u|\to+\infty$, so that the basic moment condition $M_0(J_{N,\alpha})<+\infty$ is certainly satisfied for every $N\in\mathbb{N}$. Therefore, considering Durrmeyer sampling type operators of the form $\mathcal{D}_w^{J_{N,\alpha},J_{N,\alpha}}$, the case $p=+\infty$ of Corollary \ref{corex} holds if $N>3/2$, while the case $p>1$ is valid for $N>p+1/2$. Since the case $p=1$ fails, it can be recovered by considering, for instance, operators of the form $\mathcal{D}_w^{J_{N,\alpha},\sigma_n}$.
\begin{remark}
As already observed in the introduction, taking as $\psi$ the characteristic function of the interval $[0,1]$, denoted by $\chi_{[0,1]}$, we reduce to the remarkable case of Kantorovich sampling type operators (see \cite{angamuthu2020direct,CCV23}) of the form \begin{equation*}\label{KSO}
 (\mathcal{K}_w^{\varphi}f)(x):=\sum_{k\in\mathbb{Z}}\varphi(wx-k)w\int_{\frac{k}{w}}^{\frac{k+1}{w}}f\left(u\right)du=\left(\mathcal{D}_w^{\varphi,\chi_{[0,1]}}f\right)(x),\qquad x\in\mathbb{R},
\end{equation*}
$w>0$, where $\varphi$ satisfies the basic properties of a discrete kernel. In this case, we point out that $\widetilde{m}_1(\chi_{[0,1]})=\frac{1}{2}\ne0$, so that the vanishing moment condition (ii) is not satisfied for $r=2$. In general, it is easy to see that $\widetilde{m}_\nu(\chi_{[0,1]})=\frac{1}{\nu+1}\ne0$, for every $\nu\in\mathbb{N}$. Hence, in order to improve the order of approximation up to 2, we should replace $\chi_{[0,1]}$ with a symmetric characteristic functions of the form $\dfrac{\chi_{[-a,a]}}{2a}$, with $a>0$, obtaining again operators of Kantorovich type. In fact, in this case all the continuous moments of order odd vanish.
\end{remark}
Now, we want to increase the order of approximation up to an integer $r>2$. To this aim, we have to consider suitable linear combinations of shifted kernels (see, e.g., \cite{MultiDurr}). 
For instance, we may define the following kernel (see Figure \ref{spline_traslate})
\[
\tau(u):=\frac{1}{8}\left\{47\sigma_3(u-2)-62\sigma_3(u-3)+23\sigma_3(u-4)\right\},\qquad u\in\mathbb{R},
\]
that is precisely a linear combination of shifted B-spline kernels of order $n=3$: here the coefficients are established in a such a way that the maximum value of $r\in\mathbb{N}$ for which the moments $m_\nu(\tau, x)$, $x\in\mathbb{R}$, and $\widetilde{m}_\nu(\tau)$ vanish, for every $\nu=1,\dots,r-1$, is exactly equal to 3, as a consequence of PSF (for details, see \cite{BS93}). 

\begin{remark}
Besides having a higher order of approximation than $\mathcal{D}_w^{\sigma_n,\sigma_n}$, operators based upon linear combinations of B-splines $\mathcal{D}_w^{\tau,\tau}$ also enjoy the so-called \textit{prediction property}. Being the support $[1/2,11/2]\subset (0,+\infty)$, for any $x\in\mathbb{R}$ be fixed, the series in (\ref{Durrmeyeroperator}) reduces to consider only those $k\in\mathbb{Z}$ for which $1/2\le wx-k\le 11/2$. This means that one only needs 5 samples $k/w$ taken in the past, i.e., $k/w<x$, to predict the value of the signal $f$ at the instant $x$, namely $f(x)$.
\end{remark}
We point out that the rate of convergence can be further enhanced by considering other linear combinations of B-spline $\sigma_n$, when $n>3$ (see, e.g., \cite{12}).


\section*{Acknowledgments}

{\small The authors are members of the Gruppo Nazionale per l'Analisi Matematica, la Probabilit\`a e le loro Applicazioni (GNAMPA) of the Istituto Nazionale di Alta Matematica (INdAM), of the network RITA (Research ITalian network on Approximation), and of the UMI (Unione Matematica Italiana) group T.A.A. (Teoria dell'Approssimazione e Applicazioni). 
}

\section*{Funding}

{\small The authors D. Costarelli and M. Piconi have been supported within the project: PRIN 2022 PNRR: ``RETINA: REmote sensing daTa INversion with multivariate functional modeling for essential climAte variables characterization", funded by the European Union under the Italian National Recovery and Resilience Plan (NRRP) of NextGenerationEU, under the Italian MUR (Project Code: P20229SH29, CUP: J53D23015950001). 
\\
The author G. Vinti has been supported within the project: PRIN 2022: “EXPANSION - EXPlainable AI through high eNergy physicS for medical Imaging in ONcology" funded by the European Union under the Italian National Recovery and Resilience Plan (NRRP) of NextGenerationEU, under the MUR (Project Code: 2022Z3BT9F, CUP: J53D23002530006) and by University Research Fund 2024 of the University of Perugia.
}

\section*{Conflicts of interests/Competing interests}
{\small The authors declare that they have not conflict of interest and/or competing interest.}

\section*{Availability of data and material and Code availability}

{ \small Not applicable.}


\begin{thebibliography}{99}

\bibitem{acar2024characterization}
T. Acar and B. R. Draganov, 
\textit{A characterization of the rate of the simultaneous approximation by generalized sampling operators and their Kantorovich modification}, 
Journal of Mathematical Analysis and Applications, 
\textbf{530} (2024), no. 2, 127740.

\bibitem{ABR}
A. M. Acu, I. C. Buscu, and I. Rasa, 
\textit{Generalized Kantorovich modifications of positive linear operators}, 
Mathematical Foundations of Computing, 
\textbf{6} (2023), no. 1.

\bibitem{AHR}
A. M. Acu, L. Hodis, and I. Rasa, 
\textit{Multivariate weighted Kantorovich operators}, 
Mathematical Foundations of Computing, 
\textbf{3} (2020), no. 2, 117--124.

\bibitem{ACP2012}
J. M. Aldaz, L. Colzani, and J. P{\'e}rez L{\'a}zaro, 
\textit{Optimal bounds on the modulus of continuity of the uncentered Hardy--Littlewood maximal function}, 
Journal of Geometric Analysis, 
\textbf{22} (2012), no. 1, 132--167.

\bibitem{angamuthu2020direct}
S. K. Angamuthu and S. Bajpeyi, 
\textit{Direct and inverse results for Kantorovich type exponential sampling series}, 
Results in Mathematics, 
\textbf{75} (2020), no. 3, 119.


\bibitem{ACCSV}
L. Angeloni, N. {\c{C}}etin, D. Costarelli, A. R. Sambucini, and G. Vinti, 
\textit{Multivariate sampling Kantorovich operators: quantitative estimates in Orlicz spaces}, 
Constructive Mathematical Analysis, 
\textbf{4} (2021), no. 2, 229--241.

\bibitem{bajpeyi2022approximation}
S. Bajpeyi, A. S. Kumar, and I. Mantellini, 
\textit{Approximation by Durrmeyer type exponential sampling operators}, 
Numerical Functional Analysis and Optimization, 
\textbf{43} (2022), no. 1, 16--34.

\bibitem{2}
C. Bardaro, L. Faina, and I. Mantellini, 
\textit{Quantitative Voronovskaja formulae for generalized Durrmeyer sampling type series}, 
Mathematische Nachrichten, 
\textbf{289} (2016), no. 14-15, 1702--1720.

\bibitem{12}
C. Bardaro and I. Mantellini, 
\textit{Asymptotic expansion of generalized Durrmeyer sampling type series}, 
Jaen Journal on Approximation, 
\textbf{6} (2014), no. 2, 143--165.


\bibitem{Becker79}
M. Becker, 
\textit{An elementary proof of the inverse theorem for Bernstein polynomials}, 
Aequationes Mathematicae, 
\textbf{17} (1978), no. 1, 389--390.

\bibitem{BeckerNessel78}
M. Becker and R. J. Nessel, 
\textit{An elementary approach to inverse approximation theorems}, 
Journal of Approximation Theory, 
\textbf{23} (1978), no. 2, 99--103.

\bibitem{Berdysheva}
E. E. Berdysheva, 
\textit{Uniform convergence of Bernstein--Durrmeyer operators with respect to arbitrary measure}, 
Journal of Mathematical Analysis and Applications, 
\textbf{394} (2012), no. 1, 324--336.

\bibitem{Berens72}
H. Berens, G. G. Lorentz, and R. E. MacKenzie, 
\textit{Inverse theorems for Bernstein polynomials}, 
Indiana University Mathematics Journal, 
\textbf{21} (1972), no. 8, 693--708.

\bibitem{BS22}
A. Boccuto and A. R. Sambucini, 
\textit{Some applications of modular convergence in vector lattice setting}, 
Sampling Theory, Signal Processing, and Data Analysis, 
\textbf{20} (2022), no. 2, 12.

\bibitem{BS23}
A. Boccuto and A. R. Sambucini, 
\textit{Abstract integration with respect to measures and applications to modular convergence in vector lattice setting}, 
Results in Mathematics, 
\textbf{78} (2023), no. 1, 4.

\bibitem{butzer2009voronovskaya}
P. L. Butzer and H. Karsli, 
\textit{Voronovskaya-type theorems for derivatives of the Bernstein-Chlodovsky polynomials and the Szász-Mirakyan operator}, 
Commentationes Mathematicae, 
\textbf{49} (2009), no. 1.

\bibitem{FAA71}
P. L. Butzer and R. J. Nessel, 
\textit{Fourier analysis and approximation, Vol. 1}, 
Reviews in Group Representation Theory, Part A (Pure and Applied Mathematics Series, Vol. 7), 
Marcel Dekker Inc., New York, NY, USA, 1971.

\bibitem{17}
P. L. Butzer, S. Ries, and R. L. Stens, 
\textit{Approximation of continuous and discontinuous functions by generalized sampling series}, 
Journal of Approximation Theory, 
\textbf{50} (1987), no. 1, 25--39.


\bibitem{BS93}
P. L. Butzer and R. L. Stens, 
\textit{Linear prediction by samples from the past}, 
in \textit{Advanced Topics in Shannon Sampling and Interpolation Theory}, 
Springer, 1993, pp. 157--183.

\bibitem{ITSF08}
P. L. Butzer and R. L. Stens, 
\textit{Reconstruction of signals in $L^p(\mathbb{R})$-space by generalized sampling series based on linear combinations of B-splines}, 
Integral Transforms and Special Functions, 
\textbf{19} (2008), no. 1, 35--58.

\bibitem{CCV23}
M. Cantarini, D. Costarelli, and G. Vinti, 
\textit{Approximation results in Sobolev and fractional Sobolev spaces by sampling Kantorovich operators}, 
Fractional Calculus and Applied Analysis, 
\textbf{26} (2023), no. 6, 2493--2521.

\bibitem{CMGR2014}
D. C{\'a}rdenas-Morales, P. Garrancho, and I. Ra{\c{s}}a, 
\textit{Approximation properties of Bernstein--Durrmeyer type operators}, 
Applied Mathematics and Computation, 
\textbf{232} (2014), 1--8.

\bibitem{MultiDurr}
D. Costarelli, M. Piconi, and G. Vinti, 
\textit{The multivariate Durrmeyer-sampling type operators in functional spaces}, 
Dolomites Research Notes on Approximation, 
\textbf{15} (2022), no. DRNA Volume 15.5, 128--144.


\bibitem{DurrConv}
D. Costarelli, M. Piconi, and G. Vinti, 
\textit{On the convergence properties of sampling Durrmeyer-type operators in Orlicz spaces}, 
Mathematische Nachrichten, 
\textbf{296} (2023), no. 2, 588--609.

\bibitem{EstimatesDurr}
D. Costarelli, M. Piconi, and G. Vinti, 
\textit{Quantitative estimates for Durrmeyer-sampling series in Orlicz spaces}, 
Sampling Theory, Signal Processing, and Data Analysis, 
\textbf{21} (2023), no. 1, 3.

\bibitem{DurrRegularization}
D. Costarelli, M. Piconi, and G. Vinti, 
\textit{On the regularization by Durrmeyer-Sampling type operators in ${L}^p$-spaces via a distributional approach}, 
Journal of Fourier Analysis and Applications, 
\textbf{31} (2025), no. 11, doi:10.1007/s00041-024-10121-y.

\bibitem{2022JFAA}
D. Costarelli and G. Vinti, 
\textit{Approximation properties of the sampling Kantorovich operators: regularization, saturation, inverse results and Favard classes in ${L}^p$-spaces}, 
Journal of Fourier Analysis and Applications, 
\textbf{28} (2022), no. 3, 49.

\bibitem{Derriennic}
M. M. Derriennic, 
\textit{Sur l'approximation de fonctions int{\'e}grables sur [0, 1] par des polyn{\^o}mes de Bernstein modifi{\'e}s}, 
Journal of Approximation Theory, 
\textbf{31} (1981), no. 4, 325--343.

\bibitem{DVL}
R. A. DeVore and G. G. Lorentz, 
\textit{Constructive approximation}, 
Springer Science \& Business Media, 
\textbf{303} (1993).

\bibitem{draganov2024characterization}
B. R. Draganov, 
\textit{A characterization of the rate of approximation of Kantorovich sampling operators in variable exponent Lebesgue spaces}, 
Revista de la Real Academia de Ciencias Exactas, Físicas y Naturales. Serie A. Matemáticas, 
\textbf{118} (2024), no. 2, 1--26.

\bibitem{Gonska}
H. Gonska, D. Kacs{\'o}, and I. Ra{\c{s}}a, 
\textit{The genuine Bernstein--Durrmeyer operators revisited}, 
Results in Mathematics, 
\textbf{62} (2012), 295--310.

\bibitem{HR}
M. Heilmann and I. Ra{\c{s}}a, 
\textit{A nice representation for a link between Baskakov-and Sz{\'a}sz--Mirakjan--Durrmeyer operators and their Kantorovich variants}, 
Results in Mathematics, 
\textbf{74} (2019), 1--12.

\bibitem{Heilmann}
M. Heilmann and M. Wagner, 
\textit{The Genuine Bernstein--Durrmeyer Operators and Quasi-Interpolants}, 
Results in Mathematics, 
\textbf{62} (2012), no. 3, 319--335.

\bibitem{J72}
H. Johnen, 
\textit{Inequalities connected with the moduli of smoothness}, 
Matemati{\v{c}}ki Vesnik, 
\textbf{9} (1972), no. 56, 289--305.

\bibitem{karsli2021approximation}
H. Karsli, 
\textit{On approximation to discrete q-derivatives of functions via q-Bernstein-Schurer operators}, 
Mathematical Foundations of Computing, 
\textbf{4} (2021), no. 1, 15.


\bibitem{Lorentz53}
G. G. Lorentz, 
\textit{Bernstein polynomials}, 
American Mathematical Society, (2012).

\bibitem{Lorentz63}
G. G. Lorentz, 
\textit{The degree of approximation by polynomials with positive coefficients}, 
Mathematische Annalen, 
\textbf{151} (1963), no. 3, 239--251.

\bibitem{1981jws}
L. Schumaker, 
\textit{Spline functions: basic theory}, 
Cambridge University Press, (2007).

\bibitem{Stein1970}
E. M. Stein, 
\textit{Singular integrals and differentiability properties of functions}, 
Princeton University Press, (1970).










\end{thebibliography}
\end{document}